\documentclass[11pt]{article}
\usepackage{bbm}
\usepackage{epsfig}
\usepackage{amssymb,amsmath,amscd, mdwmath}
\usepackage{graphicx,cite,url}

\begin{document}
\title{\bf Westervelt Equation Simulation on Manifold using  DEC}

\author{
Zheng Xie$^1$\thanks{E-mail: lenozhengxie@yahoo.com.cn Tel./fax: +86
0739 5316081}~~~~Yujie Ma$^2$\thanks{ E-mail: yjma@mmrc.iss.ac.cn
This work is partially supported by CPSFFP (No. 20090460102) and
NNSFC (No. 10871170) }
\\{\small
$1.$ Center of Mathematical Sciences Zhejiang University
(310027),China}
\\ {\small $2.$ Key Laboratory of Mathematics Mechanization,}
\\ {\small  Chinese Academy of Sciences,  (100090), China}}

\date{}

\maketitle

\begin{abstract} The Westervelt
equation is a model for the propagation of finite amplitude
ultrasound.  The method  of discrete exterior calculus can be used
to solve this equation numerically. A significant advantage of this
method is that it  can be used to find numerical solutions in the
discrete space manifold and the time, and therefore is a generation
of finite difference time domain method. This algorithm has been
implemented in C++.
\end{abstract}

\vskip 0.2cm \noindent {\bf Keywords: }  Westervelt equation,
Laplacian operator, Discrete exterior calculus, Manifold, Numerical
simulation.

\vskip 0.2cm \noindent {\bf MSC(2010):}  35J05, 65M12, 65M08, 53A25.
\section{Introduction}

The nonlinear effects are important in medical ultrasound and also
plays a role in diagnostic program. Some of diagnostic ultrasound
instruments  have implemented harmonic imaging into their devices by
receiving the harmonics in the reflected ultrasound caused by
nonlinear distortion of the signal propagating through the
biological tissue. Meanwhile, we should pay attention to observing
hygienic limits of ultrasound exposition to avoid heating the
unwanted tissue. The ability to predict the effects of nonlinear
ultrasound propagation therefore becomes important. Numerical
simulations are currently the best means of making predictions of
nonlinear ultrasound propagation. Two  of the numerical techniques
suitable for providing  solutions to the propagation problem are the
finite difference time domain(FDTD) method\cite{Cizek,yee,bossavit1}
and the discrete exterior calculus(DEC)\cite{meyer,hyman,hiptmair,
arnold,desbrun, whitney,leok,stern,xie-ye-ma}.

The Westervelt equation is a model for the propagation of finite
amplitude ultrasound, deriving from the equation of fluid motion by
keeping up to quadratic order terms\cite{Taraldsen,Hamilton}. In
this paper,       an conditional stable
 scheme for this equation
  using the techniques in DEC is proposed.
   A significant advantage of this scheme is that it can
be used to find  numerical solutions on the discrete space manifold
and the time,  therefore  is a generation of FDTD. In order to
refine the space mesh to improve the accuracy of the solution given
by the this scheme, the amount of work involved increases rapidly,
since the length of time step should  also be reduced.  Hence,  two
unconditional stable schemes are also proposed, namely the implicit
and semi-implicit DEC schemes for Westervelt equation on manifold.

\section{DEC schemes for Westervelt equation}

\subsection*{Discrete Laplace operator}

The $2$D or $3$D   space manifold can be approximated by triangles
or tetrahedrons,
  and the time by line segments.   Suppose each simplex contains its circumcenter. The
circumcentric dual cell $D(\sigma_0)$ of simplex $\sigma_0$ is
 $$ D(\sigma_0):=\bigcup_{\sigma_0\in \sigma_1\in\cdots \in\sigma_r}
 \mathrm{Int}(c(\sigma_0)c(\sigma_1)\cdots c(\sigma_r) ),$$
where   $\sigma_i$ is all the simplices which contains
$\sigma_0$,..., $\sigma_{i-1}$, and $c(\sigma_i)$ is the
circumcenter of $\sigma_i$. A discrete differential $k$-form, $k \in
\mathbb{Z}$, is the evaluation   of the differential $k$-form on all
$k$-simplices. Dual forms, i.e., forms  evaluated on the dual cell.
In DEC, the exterior derivative $d$ is approximated as the transpose
of the incidence matrix of $k$-cells on $k+1$-cells,
  the approximated Hodge Star $\ast$ scales the cells by
  the volumes of the corresponding dual and
primal cells, and the Laplace operator is approximated as
  $$\Delta\approx  \ast^{-1}d^{T}\ast + d^{T}\ast d.$$
Take Fig.1  as an example for a part of 2D mesh, in which $0$,...,
$C$ are vertices, $1$, $2$, $3$ are the circumcenters of triangles,
$a$, $b$, $c$ are the circumcenters of edges. Denote $l_{ij}$ as the
length of line segment $(i,j)$ and $A_{ijkl}$ as the area of
quadrangle $(i,j,k,l)$.
 $$
\begin{minipage}{0.99\textwidth}
\begin{center} \includegraphics[scale=0.45]{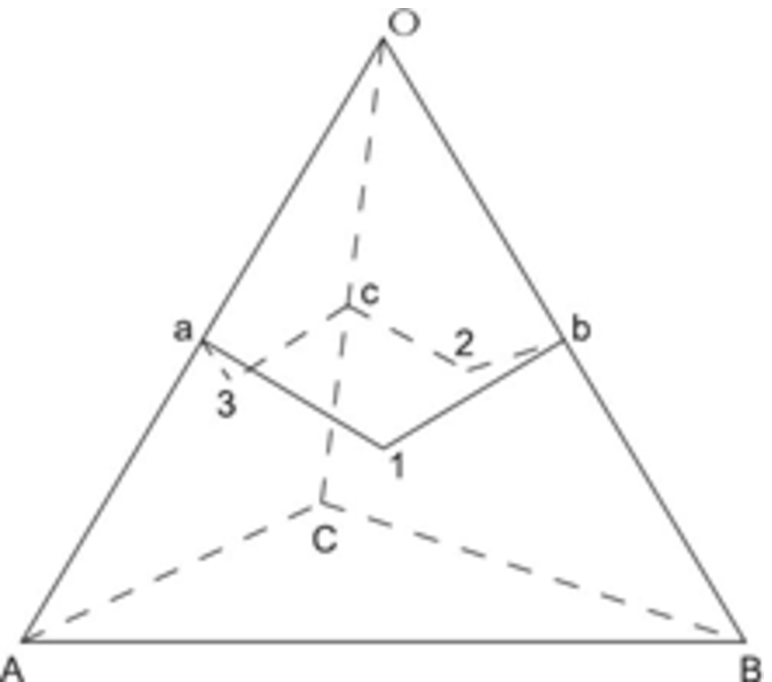}
\end{center}
\centering{ Figure.1. A part of 2D mesh}
\end{minipage}
$$
Define $$l_{12}:=l_{1b}+l_{2b},~l_{23}:=l_{2c}+l_{3c},~
l_{31}:=l_{3a}+l_{1a},$$  and $$ P_{123}:=
A_{01ab}+A_{02bc}+A_{03ac} .$$ In Fig.1, the discrete Laplace
operator acting on $p$ at vertex $0$
  is
$$\Delta p_0 \approx \frac{ 1} {P_{123 }} \left(\frac{l_{13}}{l_{A0}}
(p_A-p_0)+\frac{l_{12}}{l_{B0}}(p_B-p_0)+\frac{l_{23}}{l_{C0}}(p_C-p_0)
\right) . \eqno{(1)}
$$
$$
\begin{minipage}{0.99\textwidth}
\begin{center} \includegraphics[scale=0.45]{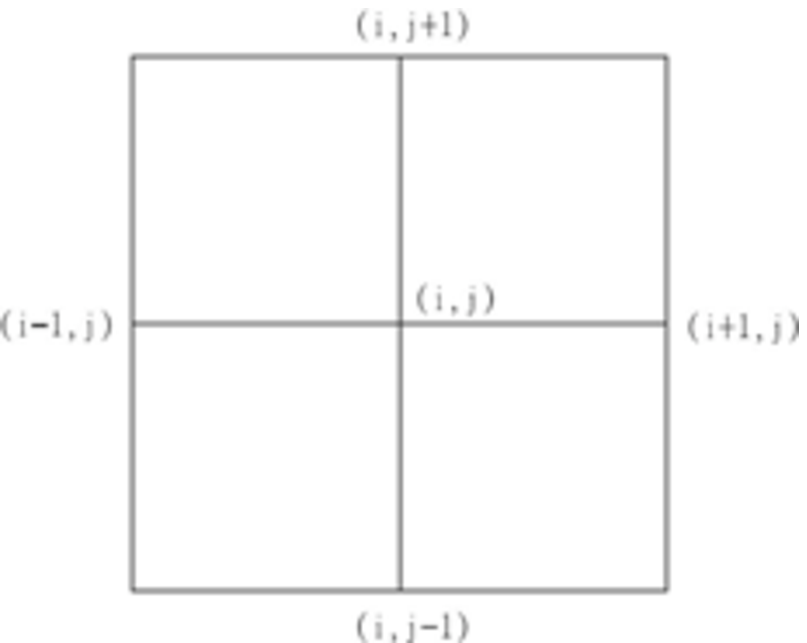}
\end{center}
\centering{ Figure.2. A part of 2D rectangular mesh}
\end{minipage}
$$   In Fig.2,  scheme (1)
reduces to
$$
\Delta p_{i,j}\approx\dfrac{ p_{i,j+1} +p_{i,j-1}
 +p_{i+1,j }+p_{i-1,j}-4p_{i,j } }{(\Delta
 s)^2}, $$where  $\Delta s$ is the uniform space step length.

\subsection*{Discretization for Westervelt equation}

The Weservelt equation can be written in the following form:
$$\Delta p-\frac{1}{c^2_0}\frac{\partial^2p}{\partial
t^2}+\frac{\delta}{c^4_0} \frac{\partial^3p}{\partial
t^3}+\frac{\beta}{\rho_0c^4_0} \frac{\partial^2p^2}{\partial t^2}=0,
\eqno{(2)}$$
 where $p$ is the acoustic pressure,
$\rho_0$ is the are the ambient density, and $c_0$ is the ultrasound
speed, $\delta$ is the diffusivity of ultrasound, $\beta$ is the
coefficient of nonlinearity.

For some situations, a source having azimuthal symmetry about its
axis is considered. In this case,  2D triangular discrete manifold
as the space is only need to be considered. And the Eq.(2) on 3D
space manifold and the time  can also be solved numerically, using a
similar approach.

The temporal partial
 derivatives with discrete differences, which can be obtained from
 Taylor series expansions about each node of the computational mesh.
Temporal derivatives present in Eq.(2) can be calculated  as
follows:
$$\begin{array}{lll}\dfrac{\partial^2p^n}{\partial t^2}\approx \dfrac{1}{(\Delta t)^2}(p^{n+1}-2p^n+p^{n-1})\\
\dfrac{\partial^3p^n}{\partial t^3}\approx \dfrac{1}{(\Delta
t)^3}(p^{n}-3p^{n-1}+3p^{n-2}-p^{n-3})
\\ \dfrac{\partial^2(p^2)^n}{\partial
t^2}\approx\dfrac{1}{(\Delta
t)^2}((p^2)^{n}-2(p^2)^{n-1}+(p^2)^{n-2}), \end{array}$$ where
$\Delta t$ is the time step length, $n$ is
 the time coordinate. Let
$$\begin{array}{lll}L p^{n-1}_0&:=&\frac{1}{(c_0\Delta
t)^2}(p^{n}_0-2p^{n-1}_0+p^{n-2}_0) -\frac{\delta}{c^4_0(\Delta
t)^3}(p^{n-1}_0-3p^{n-2}_0+3p^{n-3}_0 \\&&-p^{n-4}_0)-
\frac{\beta}{\rho_0c^4_0(\Delta
t)^2}((p^2_0)^{n-1}-2(p^2_0)^{n-2}+(p^2_0)^{n-3})\end{array}$$ The
explicit DEC scheme for Eq.(2) is
$$  \frac{  \frac{l_{13}}{l_{A0}}
(p^{n-1}_A-p^{n-1}_0)+\frac{l_{12}}{l_{B0}}(p^{n-1}_B-p^{n-1}_0)+\frac{l_{23}}{l_{C0}}(p^{n-1}_C-p^{n-1}_0)
 } {P_{123 }}=L p^{n-1}_0. \eqno{(3)}$$

If refining the space mesh, the amount of work involved increases
rapidly, since   the length of time step   should  satisfy the
stable condition. Now,  two unconditional stable schemes are
proposed, namely the implicit DEC scheme (4) and semi-implicit DEC
scheme (5).
$$ \frac{ 1} {P_{123 }} \left(\frac{l_{13}}{l_{A0}}
(p^{n}_A-p^{n}_0)+\frac{l_{12}}{l_{B0}}(p^{n}_B-p^{n}_0)+\frac{l_{23}}{l_{C0}}(p^{n}_C-p^{n}_0)
\right)=L p^{n-1}_0 \eqno{(4)}$$
$$ \frac{ 1} {P_{123 }}
\left(\frac{l_{13}}{l_{A0}}
(p^{n-1}_A-p^{n}_0)+\frac{l_{12}}{l_{B0}}(p^{n-1}_B-p^{n}_0)+\frac{l_{23}}{l_{C0}}(p^{n-1}_C-p^{n}_0)
\right)=L p^{n-1}_0 \eqno{(5)}
$$

The higher order accuracy scheme for temporal derivatives can also
be used in schemes(3-5).

\section{Stability, convergence  and accuracy}

The first two terms in (2) describe linear  wave propagation at the
small-signal sound speed. The third term describes the loss due to
he viscosity and thermal conduction of the media. The fourth term
describes the nonlinear distortion of the traveling wave due to
amplitude effects. Since the coefficients of third and fourth terms
in Eq.(2)  are small compared with the coefficient of second term,
   the effect of $\frac{\partial^3p}{\partial t^3}$ and
$\frac{\partial^2p^2}{\partial t^2}$ on the stability of the scheme
(3) can be ignored. Therefore,  the stable condition for   the wave
equation
$$\Delta p-\frac{1}{c^2_0}\frac{\partial^2p}{\partial
t^2}=0   $$ is    need to be satisfied. The wave equation can be
discretized as
$$\begin{array}{lll}\mathrm{Right}(1)^{n-1}&=&\dfrac{1}{(c_0\Delta
t)^2}(p^{n}_0-2p^{n-1}_0+p^{n-2}_0)
\end{array}\eqno{(6)}$$
The stable condition condition for scheme (6)\cite{xie-ye-ma} is
$$ c_0\Delta t \leq \mathrm{min}_{\forall 0\in \mathrm{vertices}}
\sqrt{\dfrac{2}{\dfrac{1}{P_{123}}{\left(\dfrac{l_{13} }{l_{A0} }
+\dfrac{l_{12} }{l_{B0} } +\dfrac{l_{23} }{l_{C0}
}\right)}}}.\eqno{(7)}$$
 In Fig.2, the inequality (7) reduces
to
$$c_0\Delta t \leq \dfrac{\sqrt{2}}{2} \Delta s,$$
which is the stable condition for scheme (6) on square grid. The
analysis of   unconditional stability for schemes (4) and (5) can
also be proved as the wave equations.

By the definition of truncation error, the exact solutions
  of Westervelt equation  satisfy the same relation as schemes (3-5)
  except for an additional term $O (\Delta t)^2  $.
   This expresses the consistency, and so  the convergence for  schemes (3-5) by  Lax equivalence theorem
   (consistency $+$ stability $=$ convergence).
 The  derivative in Eq.(2) is approximated by first order
difference in schemes (3-5). Equivalently, $p$ is approximated by
linear interpolation function. This is to say  that schemes (3-5)
have first order spacial and temporal accuracy, according to the
  definition about accuracy of finite volume method.
\section{Algorithm Implementation}

The DEC schemes (3-5) for Westervelt equation was implemented in
C++. The implementation   consists of the following steps:
\begin{itemize}
  \item [1.]Set the simulation parameters. These
are the dimensions of the computational mesh and the size of the
time step, etc.;

  \item [2.] Set the propagating media parameters.

 \item [3.] Initialize the mesh indexes.

  \item [4.]Assign current transmitted signal.

  \item [5.]Compute the value of all spatial nodes  and temporarily store the
  result  in the circular buffer for further computation.

  \item [6.]Visualize the currently computed grid of spatial nodes.

   \item [7.]Repeat the whole process  from the step 4, until reach the desired total number
   of iterations.
\end{itemize}

In the common practice, not every simulation step needs to be
visualized, especially when the time step size is too small. In all
the following examples, scheme (3) is  used, and  the parameters are
$\delta=0.01$ $\beta=1$ $\rho=10000$. The Gaussian envelope $
\frac{1}{\sqrt{2\pi}}\exp(-{\frac{t^2}{2}})\cos(t)$ is used as a
source signal.

The example in Fig.3 shows the nonlinear propagation effects and the
diffusivity in the propagation of the Gaussian  envelope.
$$
\begin{minipage}{0.99\textwidth}
\begin{center} \includegraphics[scale=0.45]{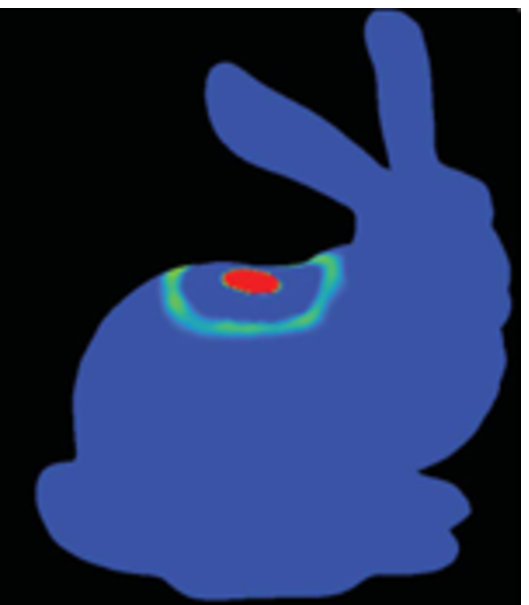}
 \includegraphics[scale=0.45]{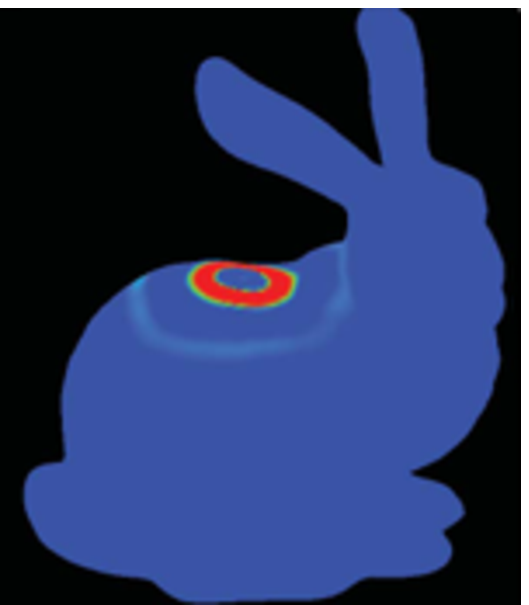}
 \includegraphics[scale=0.45]{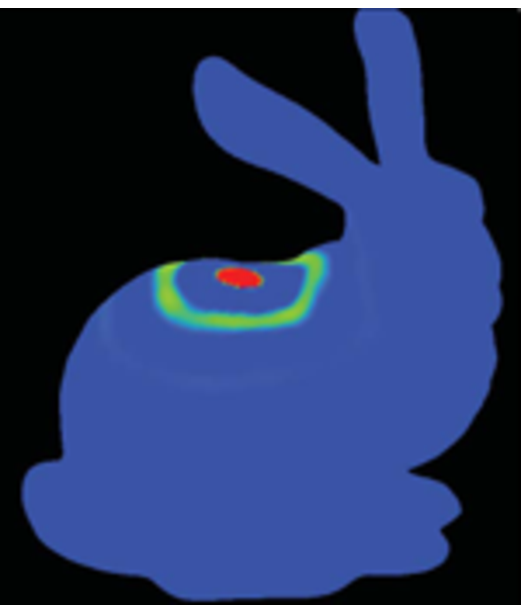}
\end{center}
\centering{ Figure.3. The propagation of  Gaussian  envelope  on a
rabbit}
\end{minipage}
$$

The example in Fig.4 shows the propagation of Gaussian  envelope at
a boundary of two kinds of media with different wave speeds $340$m/s
and $3400$m/s.

$$
\begin{minipage}{0.99\textwidth}
\begin{center} \includegraphics[scale=0.45]{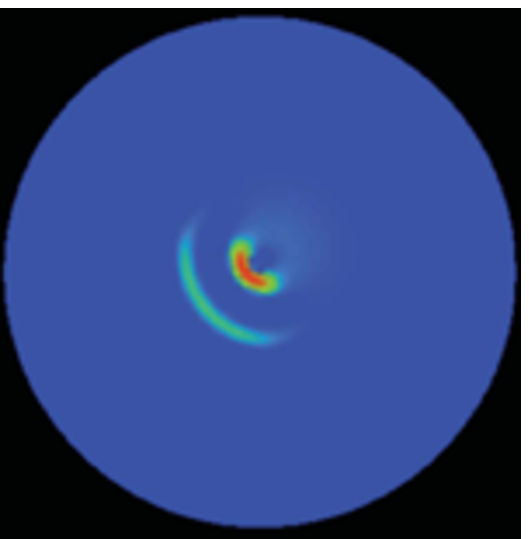}
 \includegraphics[scale=0.45]{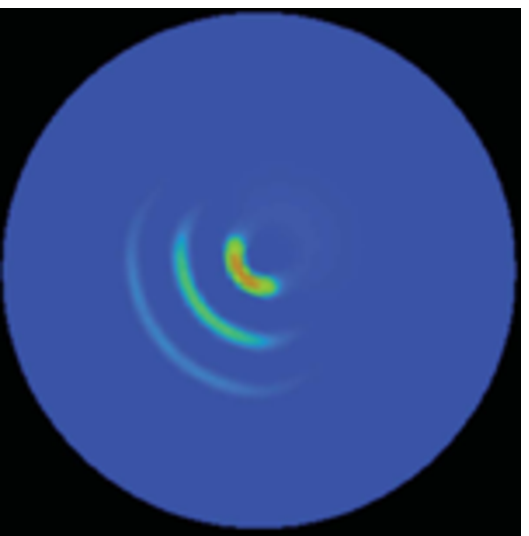}
 \includegraphics[scale=0.45]{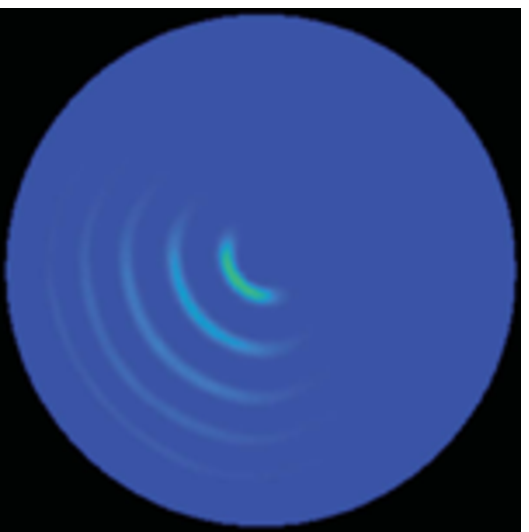}
\end{center}
\centering{Figure.4.  The propagation of   Gaussian  envelope  on a
sphere with two kinds of media }
\end{minipage}
$$


\begin{thebibliography}{99}

\bibitem {yee} K.S. Yee,
Numerical solution of inital boundary value problems involving
Maxwell's equations in isotropic media, IEEE Trans. Ant. Prop.,
14(3), 302-307, (1966).





\bibitem {bossavit1} A. Bossavit, L. Kettunen,  Yee-like schemes on a
tetrahedral mesh, with diagonal lumping, Int. J. Numer. Modell.,
12(1-2), 129-142, (1999).

 \bibitem {Cizek} M. Cizek,  J. Rozman,   Acoustic Wave Equation Simulation Using
FDTD, Radioelektronika,  17th International Conference, (2007).




 \bibitem {whitney} H. Whitney, Geometric integration theory,
Princeton University Press, Princeton, (1957).



 \bibitem {arnold} D.N. Arnold, R.S. Falk, R. Winther, Finite element exterior calculus,
homological techniques, and applications, Acta Numer., 15, 1-155,
(2006).




 \bibitem {meyer}M. Meyer, M. Desbrun, P. Schr\"oder,  A.H. Barr,
Discrete differential geometry operators for triangulated
2-manifolds. In InternationalWorkshop on Visualization and
Mathematics, VisMath, (2002).

 \bibitem {desbrun} M. Desbrun, A.N. Hirani, M.
Leok, J. E. Marsden, Discrete exterior calculus, arXiv:
math.DG/0508341.

 \bibitem {hyman}J. M. Hyman, M. Shashkov, Natural discretizations for the
divergence, gradient, and curl on logically rectangular grids,
Comput. Math. Appl., 33(4):81-104, (1997).


 \bibitem {hiptmair}R. Hiptmair, Discrete Hodge operators, Numer. Math., 90(2), 265-289,
(2001).




\bibitem {leok} M. Leok, Foundations of computational geometric mechanics.
Ph.D. thesis, California Institute of Technology, (2004).

\bibitem {stern} A.
Stern, et al, Computational Electromagnetism with Variational
Integrators and Discrete Differential Forms, arXiv:0707.4470.

\bibitem{xie-ye-ma} Z. Xie, Y.J. Ma,   Computation of Maxwell's equations on Manifold using DEC,
 arXiv:0908.4448.


 \bibitem {Hamilton} M.F. Hamilton, C.L. Mofrey,  Model equations, Nonlinear
Acoustics, Chap.3, Academic Press, San Diego, (1998).


 \bibitem {Taraldsen}G. Taraldsen, A generalized Westervelt equation for nonlinear
medical ultrasound, J. Acoust. Soc. Am. 109(4), 1329-1333, (2001).






\end{thebibliography}
\end{document}